\documentclass[11pt]{article}
\usepackage{geometry}                % See geometry.pdf to learn the layout options. There are lots.
\usepackage{graphicx}
\usepackage{amssymb}
\usepackage{epstopdf}
\title{A low memory, highly concurrent multigrid algorithm}
\author{Mark F. Adams\thanks{Applied Physics and Applied Mathematics Department, Columbia University, \texttt{mark.adams@columbia.edu}}}
\date{}                                           % Activate to display a given date or no date

\newcommand{\Order}[1]{\ensuremath{\mathcal{O}(#1)}}    % big O notation

\begin{document}

\maketitle

\begin{abstract} 
We examine what is an efficient and scalable nonlinear solver, with low work and memory complexity, for many classes of discretized partial differential equations (PDEs) -- matrix-free Full multigrid (FMG) with a Full Approximation Storage (FAS) -- in the context of current trends in computer architectures.
Brandt proposed an extremely low memory FMG-FAS algorithm in the 1970s that has several attractive properties for reducing costs on modern -- memory centric -- machines and has not been developed.
This method, \textit{segmental refinement} (SR), has very low memory requirements because the finest grids need not be held in memory at any one time but can be ``swept" through, computing coarse grid correction and any quantities of interest, allowing for orders of magnitude reduction in memory usage.
This algorithm has two useful ideas for effectively exploiting future architectures: improved data locality and reuse via ``vertical" processing of the multigrid algorithms and the method of $\tau$-corrections, which allows for not storing the fine grid(s).
This report develops a parallel generalization of the original sweeping technique and explores algorithmic details with the 1D model problem. 
We show that FMG-FAS-SR can work as originally predicted, solving systems accurately enough to maintain the convergence rate of the discretization with one FMG iteration, and that the parallel algorithm provides a natural approach to fully exploiting the available parallelism of FMG.
The parallel algorithm is naturally expressed in asynchronous data-driven programming models, which is responsive to current directions in programming models for extreme scale machines.
\end{abstract}

\section{Motivation}

Current trends in computer architecture such as non-rising, or even falling, clock rates, saturated processor architectures in terms of pipelining, etc., the continued increase in the number of transistors on a chip, requires that algorithms be highly concurrent and asynchronous.
Additionally, we have reached a point where the exponential growth in power cost, which goes along with this continued growth of extreme-scale machines, is becoming the prohibitive  cost of extreme-scale computing.
The desire to continue the exponential growth of extreme-scale PDE simulations combined with an economic need to keep the power budget of a machine down to say 25MW will tax the resources of computer engineers and may require that we develop algorithms for radically different machine models with respect to memory, energy and faults than what we have worked with in the past.
These changes, along with the continued need for mathematically scalable algorithms, as we increase the fidelity of our simulations, is leading to the need to rethink our solver algorithms for large scale PDE simulations.
In particular, the powering and moving of memory will become  more central to the cost of PDE solves and the flop counts will become less so.
This paper aims to address the root cause of these future costs -- memory -- by developing low memory and memory movement PDE solver algorithms that exploit the mathematical structure of PDEs.

%There are many levels of techniques in the system stack that will be deployed to achieve the goal of producing real science on an exa-scale machine, from hardware, firmware, systems software, and programming models, to computer science algorithms (eg, sorting), mathematical algorithms (eg. solvers), discretizations and the mathematical models themselves.
To rationally develop an algorithm, and certainly to analyze an algorithm, one needs a machine model.
Traditional complexity theory, essentially counting operations or flops, has served this purpose well for high performance computing -- it along with its extensions to parallel complexity -- has a well developed theory.
Memory complexity is also useful and to some extent serves as a proxy for memory movement complexity.
While data locality, to reduce memory traffic in the memory hierarchy, has been central to high performance computing for decades it is difficult to incorporate memory movement into complexity models directly and there is no consensus on any one approach though much work has been done in this area \cite{Aggarwal:1987am,Aggarwal:1988hx,Savage:2008bk,Ballard:2009ao,Blelloch:2010qf}.
The dearth of good cost models for future machines, whose design is an active area of research and far from well understood, leads us to look at the fundamental source of costs -- memory -- and place less emphasis on what has historically been the primary measure of costs -- flops.

Multigrid is an efficient solver method for many classes of problems; matrix-free Full multigrid (FMG), with a Full Approximation Storage (FAS) for nonlinear problems, is an very efficient algorithm for some classes of problems \cite{UTrottenberg_CWOosterlee_ASchueller_2000a}, with very low memory and work complexity.
Brandt proposed an extremely low memory FMG algorithm in the 1970s -- segmental refinement (SR)(\cite{ABrandt_1977a} \S 7.5; \cite{ABrandt_1979b} \S 2.2; \cite{Brandt-2011} \S 8.7; \cite{NDinar_1979a};\cite{ABrandt_BDiskin_1994a}) -- that does not require that all of the data be stored at any one time.
The resulting serial algorithm has a memory complexity of $\log^{D+1}{(N)}$  \cite{Brandt-2011} in $D$ dimensions.
This is done by reformulating the multigrid algorithms from the view of the coarse grid, using the method of $\tau$-corrections, so that in effect of the fine grid is stored (compressed) on the coarse grid and can be recovered with a special multigrid smoothing technique \S\ref{ssec:Kaczmarz}.
Though the algorithm was originally proposed as a low memory serial method that ``sweeps" across the grid, and was developed as a parallel algorithm by Brandt and Diskin \cite{ABrandt_BDiskin_1994a}.
%we propose a parallel variant that uses ``patches" and Jacobi (additive) smoothers \cite{Smith-96}.
This approach in effect allows for low memory to be traded for concurrency -- the coarse grids are stored explicitly and behave like standard FMG, while the finer grids do not explicitly store the entire domain at any one time.
This results in higher flop to memory ratios than traditional multigrid methods because at some level of granularity, say a uniform memory access partition, the same memory is used repeatedly for many patches and thus many flops.
Additionally, this algorithm requires that the multigrid algorithm be processed ``vertically", which \textit{maximizes data reuse and locality}, as opposed to the traditional ``horizontal" implementation approach where entire grids are processed sequentially (ie, entire grids are first smoothed, residuals are calculated, and then restricted to coarse grids, etc. \S\ref{sec:mgv}).
This algorithm is also inherently asynchronous and could be naturally expressed in an asynchronous, task oriented programming language although this is not neccessary.
This algorithm has not been developed because it requires more flops than the traditional approach and is more complex to engineer, but it has many attractive properties on more memory centric computers.
%This paper explores a point in the design space of multigrid algorithms that may be useful for the more memory centric extreme-scale machines of the future.

This paper proceeds by providing some basic multigrid background in \S\ref{sec:mg}, the segmental refinement algorithm is developed in \S\ref{sec:algo} along with a parallel FMG-FAS-SR algorithm.
We apply this algorithm to a model problem in \S\ref{sec:numr} and conclude in \S\ref{sec:conc}.

\section{Multigrid Background}
\label{sec:mg}

Multigrid is an effective method for solving systems of algebraic equations that arise from discretized PDEs.
Modern multigrid's antecedents can be traced back to Southwell in the 1930s \cite{RVSouthwell_1940a}, and Fedorenko in the early 1960s \cite{RPFedorenko_1961a}.
Brandt developed multigrid's modern form in the 1970s -- algorithms and analysis with work complexities equivalent to a few residual calculations (work units), applied to complex domains, non-constant coefficients problems and nonlinear problems \cite{ABrandt_1973a}.
A substantial body of literature, both theoretical and experimental exists that proves and demonstrates the optimality of multigrid, having \Order{n} work complexity and \Order{\log(n)} parallel work complexity or computational depth for the Laplacian \cite{UTrottenberg_CWOosterlee_ASchueller_2000a}.
Multigrid has been applied to a wide range of problems \cite{Brandt-2011,UTrottenberg_CWOosterlee_ASchueller_2000a}, starting with flow problems in the seminal paper by Brandt \cite{ABrandt_1973a}.
Multigrid as also been found to useful as a nonlinear solver -- used directly on the nonlinear system -- with demonstrated costs very similar to that of a linear multigrid solve (\cite{ABrandt_1973a,UTrottenberg_CWOosterlee_ASchueller_2000a} \S 5.3.3).

%Multigrid is often understood as an algebraic equation solver for an input fine grid, where the coarse grid accelerates convergence of a simple solver like Gauss-Seidel on the fine grid, but it is beneficial to look at the fine grid as a smoother for a coarse grid solution.
%This dual point of view is also more natural if one is working within an adaptive mesh refinement (AMR) environment.
%From the coarse grid point of view one can solve the coarse grid equations with a $\tau$-correction and compute a solution of the fine grid system on the coarse grid.
%This is the basis of the methods pursued in this paper.

\subsection{Multigrid V-cycle}
\label{sec:mgv}

Multigrid methods are motivated by the observation that a low resolution discretization of an operator can capture modes or components of the error that are expensive to compute directly on a highly resolved discretization.
More generally, any poorly locally-determined solution component has the potential to be resolved with a coarser representation.  
This process can be applied recursively with a series of coarse grids, thereby requiring that each grid resolve only the components of the error that it can solve efficiently.
This process is known as a $V-cycle$ because of the shape of the graph in the standard representation of these algorithms (see Figure \ref{fig:FMG}).
These coarse grids have fewer grid points, typically about a factor of two in each dimension, such that the total amount of work in multigrid iterations can be expressed as a geometric sum that converges to a small factor of the work on the finest mesh.  
These concepts can be applied to problems with particles/atoms or pixels as well as the traditional grid or cell variables considered here.
Multigrid provides a basic framework within which particular multigrid methods can be developed for particular problems.
This framework has proven to be an effective way to separate the near-field from the far-field contributions to the solution of say elliptic operators -- the coarse grid captures the far-field contribution and the near-field is resolved with a local process called a \textit{smoother}.

The coarse grid space can be represented algebraically as the columns of the {\it prolongation} operator $I^h_{H}$, where $h$ is the fine grid mesh spacing, $H$ is the coarse grid mesh spacing. 
The prolongation operator is used to map corrections to the solution from the coarse grid to the fine grid.
Residuals are mapped from the fine grid to the coarse grid with the {\it restriction} operator $I^H_{h}$; $I^H_{h}$ is often equal to the transpose of $I^h_{H}$.
The coarse grid matrix can be formed in one of two ways, either algebraically to form Galerkin (or variational) coarse grids ($L_{H} \leftarrow I^H_{h}L_{h}I^h_{H}$) or, by creating a new operator on each coarse grid (if an explicit coarse grid is available).

\subsection{Nonlinear multigrid}

The multigrid $V$--$cycle$ can be adapted to a nonlinear method by observing that the coarse grid residual equation can be written as 
\begin{equation} 
r_{H} = L_{H} (u_{H} ) - L_{H} ({\hat u}_{H} ) = L_{H} ({\hat u}_{H} +e_{H} ) - L_{H} ({\hat u}_{H} ),
\label{eq:cresid}
\end{equation}
where $u$ is the exact solution, ${\hat u}^H$ approximates $I^H_h{u}^h$, the full intended solution represented on the coarse grid, hence the name ``Full Approximation Scheme", and $e$ is the error.
With this, and an approximate solution on the fine grid $\tilde u_h$, the coarse grid equation can be written as
\begin{equation} 
L_{H}\left (I^{H}_h {\tilde u}_h+ e_{H}\right) = L_{H} \left(I^{H}_h{\tilde u}_h\right) + I^{H}_h\left( f_h - L_h {\tilde u}_h \right), 
\label{eq:cresid2}
\end{equation}
and is solved approximately.
After $I^{H}_h {\tilde u}_h$ is subtracted from the $I^{H}_h {\tilde u}_h+ e_{H}$ term the correction is applied to the fine grid with the standard prolongation process.
This method is called Full Approximation Scheme (or Full Approximation Storage - FAS), because the full solution is stored on each level and not just a residual correction.
See Trottenberg for more details \cite{UTrottenberg_CWOosterlee_ASchueller_2000a}.

Figure \ref{fig:mgv} shows the FAS multigrid {\it V($\nu 1$,$\nu 2$)-cycle} and uses a nonlinear smoother $u \leftarrow S(L,u,f)$.
\vskip .2in
\begin{figure}[ht!]
\vbox{ \raggedright 
$\phantom{}u=${\bf  function} $FAS(L_k,u_k,f_k)$ \\ 
$\phantom{MM}${\bf if $k > 0$} \\
$\phantom{MMMM}u_k \leftarrow S^{\nu 1}(L_k,u_k,f_k)$ \quad \qquad -- $\nu 1$ iterations of the (pre) smoother \\ 
$\phantom{MMMM}r_k \leftarrow f_k - L_ku_k$ \\ 
$\phantom{MMMM}r_{k-1}\leftarrow I^{k-1}_k(r_k)$ \qquad \quad \qquad -- restriction of residual to coarse grid\\ 
$\phantom{MMMM}u_{k-1}\leftarrow I^{k-1}_k(u_k)$ \qquad \quad \qquad -- restriction of solution to coarse grid\\ 
$\phantom{MMMM}c_{k-1}\leftarrow FAS(L_{k-1},u_{k-1},r_{k-1}+L_{k-1}u_{k-1})$  \qquad -- recursive application \\
$\phantom{MMMM}u_k \leftarrow u_k + I^k_{k-1}(c_{k-1} - u_{k-1})$ \quad \qquad -- prolongation of coarse grid correction \\ 
$\phantom{MMMM}u_k \leftarrow S^{\nu 2}(L_k,u_k,f_k)$ \qquad  \qquad  -- $\nu 2$ iterations of the (post) smoother \\
$\phantom{MM}${\bf else}\\ 
$\phantom{MMMM}u_k \leftarrow L_0^{-1}f_0$ \qquad \qquad \qquad \quad -- exact solve of coarsest grid \\ 
$\phantom{MM}${\bf return} $u_k$}
\caption{FAS Multigrid {\it $V$-$cycle$} Algorithm}
\label{fig:mgv}
\end{figure}
With $M$ coarse grids the preconditioner (solver) for $L_M u_M = f_M $ is $u = FAS(L_M,0,f_M)$

\subsection{Full Multigrid}

An important variant on the $V$--$cycle$ is the $F$--$cycle$ or related full multigrid (FMG).
The multigrid $F$--$cycle$ restricts the right hand side from the fine grid to the coarsest grid and then applies a multigrid cycle, of some sort, at each level, starting with the coarsest level and interpolating the solution to the next finest level as an initial solution for the next $V$--$cycle$.
A higher order interpolator between the level solves, $\Pi^h_H$, is needed for optimal efficiency of the FMG process but requires more data movement.
An attractive property of the $F$--$cycle$ is that for some operators it has been proven that one $F$--$cycle$ is sufficient to reduce the error to the order of the truncation error, which is often all that is required \cite{Bank-81} (\cite{UTrottenberg_CWOosterlee_ASchueller_2000a} \S 3.2.2).
Thus, the algebraic system can be solved to spatial truncation error accuracy with a work complexity of a few \textit{work units}, or residual calculations.
Note, the parallel complexity of an $F$--$cycle$ does have an extra $log(n)$ factor.

One can analyze the $F$--$cycle$ with induction where the induction hypothesis is that the algebraic error is some factor $r$ of the truncation error (which is satisfied on the coarsest grid where an accurate solver is required), and the standard assumption that the truncation error is of the form \Order{h^{p}}, and that the solver on each level (eg, one $V$--$cycle$) reduces the error by some factor $\Gamma$ (which can be proven or measured experimentally) to derive an equation that directly relates $r$ to $\Gamma$.
This allows the use of the desired ratio -- any desired ratio -- of solution to truncation error to tune the solver at each level -- see Adams for the application of these ideas to compressible resistive magnetohydrodynamics where two $V$-$cycle$ were used as the level solver \cite{Adams-10a}.

FMG starts with the coarse grid, and is more natural in an AMR context; it simply omits the initial restriction of the residual to the coarse grid.
Figure \ref{fig:fmg} shows the FMG algorithm.
\vskip .2in
\begin{figure}[ht!]
\vbox{ \raggedright 
$\phantom{}u=${\bf  function} $FMG$ \\ 
$\phantom{MM}u_0 \leftarrow FAS\left(L_0, 0, f_0\right)$ \qquad \qquad \qquad \quad -- exact solve of coarsest grid \\ 
$\phantom{MM}${\bf for k=1:M} \\
$\phantom{MMMM}{u}_{k} \leftarrow \Pi^{k}_{k-1} u_{k-1}$  \qquad \qquad \qquad \quad -- FMG prolongation \\
$\phantom{MMMM}{u}_{k} \leftarrow FAS\left(L_k,u_k,f_k\right)$  \qquad \qquad -- V-cycle\\
$\phantom{MM}${\bf return} $u_M$}
\caption{ FMG-FAS algorithm}
\label{fig:fmg}
\end{figure}

FMG-FAS multigrid is an efficient solver for some classes of problems and its application to new classes of problems is an active area of research \cite{Adams-10a}.
This method was developed in the 1970s and was attractive because of its low memory requirements: only requiring the field variables themselves, and because of its very low work complexity (as low as six work units to solve to truncation error).
After the profligate era of the 1980s to the 2000s, with large amounts of uniform access memory available, low memory complexity algorithms are attractive again as we move to memory centric cost models.
Thus, FMG-FAS is an attractive solver algorithm for the anticipated machine models for exa-scale machines.

\subsection{Segmental Refinement and PDE Compression}

Looking at FAS from a two grid point of view we can rewrite the coarse grid Eq. \ref{eq:cresid2} as
\begin{equation} 
L^H{\hat u}^H = {\hat f}^H = L^H \left({\hat I}^H_h {\tilde u}^h \right) + I^H_h r^h = L^H \left({\hat I}^H_h {\tilde u}^h \right) + I^H_h \left( f^h - L^h{\tilde u}^h \right),
\label{eq:cresid3}
\end{equation}
where ${\hat I}^H_h$ is some fine-to-coarse transfer which need not be the same as $I^H_h$ (they are in principle defined on different
spaces), and  ${\tilde u}^h$ is the current solution on the fine grid.
Having obtained an approximate solution ${\tilde u}^H$ from solving Eq. \ref{eq:cresid3} we can write the fine grid correction as
\begin{equation} 
{\tilde u}^h_{NEW} =  {\tilde u}^h +I^h_H\left ( {\tilde u}^H - {\hat I}^H_h{\tilde u}^h \right).
 \label{eq:fupdate}
\end{equation}
Or use the solution directly:
\begin{equation} 
{\tilde u}^h_{NEW} =  I^h_H {\tilde u}^H.
\label{eq:fupdate2}
\end{equation}
Generally Eq. \ref{eq:fupdate2} is not preferred because it introduces interpolation error in the full solution, and not just a correction, but it will be useful in the context of this work.

%There is potential to make FAS-FMG very memory efficient by not storing the entire fine grid field variables explicitly via the method of ``segmental refinement"  \cite{Brandt-2011}.
%Brandt proposed this method as an ultra-low memory multigrid algorithm, but it also has the desirably property of containing a large amount of parallelism.
We can look at the FMG method from the dual point of view, that is from the view of the coarse grid.
Instead of looking at the coarse grid as an accelerator to the fine grid convergence we look at a fine grid as a correction to the coarse grid problem.
Eq. \ref{eq:cresid3} can be rewritten in the form:
\begin{equation} 
L^H{\hat u}^H = f^H + \tau_h^H,
\label{eq:dual_cresid}
\end{equation}
where the $\tau$-correction is 
\begin{equation} 
\tau^H_h = L^H \left( {\hat I}^H_h {\tilde u}^h\right) - I^H_h\left(L^h {\tilde u}^h\right),
\label{eq:tau}
\end{equation}
and  $f^H=I_h^Hf^h$.
At convergence ${\hat u}^H = {\hat I}^H_h {u}^h$, hence $\tau_h^H$ is the \textit{fine-to-coarse defect correction} designed to make its solution coincide with the fine-grid solution.
This observation, along with the update of Eq. \ref{eq:fupdate2} allows for an FMG-FAS algorithm that need not store the fine grids, but can compute them locally patch-by-patch.
Brandt proposed the segmental refinement method that exploits this property by ``sweeping" through the grid refining one segment at a time (\S 8.7 \cite{Brandt-2011}).

%, relax them locally, so as to compute the $\tau$-correction, again locally, one patch at a time.
Note, with the $\tau$-correction, the coarse grid solution is equal to the fine grid solution at the coarse grid points -- this allows for the inexpensive computation of the solution with a special relaxation method in the post smoothing leg of the $V$--$cycle$ (\S\ref{ssec:Kaczmarz}).
Thus, this representation can be viewed as a compression technique that exploits the PDE and multigrid method -- \textit{PDE compression}.
%This property can also be exploited to recover from compute node failures by saving (in memory) checkpoint/restart files with other compute nodes \ref{jed}.

\section{Algorithm}
\label{sec:algo}

%\documentclass[11pt]{amsart}
%\usepackage{geometry}                % See geometry.pdf to learn the layout options. There are lots.
%\geometry{letterpaper}                   % ... or a4paper or a5paper or ... 
%\usepackage{graphicx}
%\usepackage{amssymb}
%\usepackage{epstopdf}
%\DeclareGraphicsRule{.tif}{png}{.png}{`convert #1 `dirname #1`/`basename #1 .tif`.png}
%\newcommand{\Order}[1]{\ensuremath{\mathcal{O}(#1)}}    % big O notation
%\begin{document}

%Segmental refinement exploits the fact that smoothing is a local process and that errors due to boundary conditions can be ameliorated by using sufficient buffer (halo or ghost) cells around the subdomain (``patch") of interest, with approximate (eg, low order) boundary conditions on each patch.
%This removes many data dependancies that are not strictly required in the FAS-FMG algorithm.
The duel view of FMG allows the $\tau$-corrections to be computed on subdomains and the fine grid data need not be retained in memory.
In serial this allows only small parts of the fine grids to be stored at any given time as the algorithm ``sweeps" through the grids, computing the $\tau$-corrections and restricting them to the coarse grid.
This algorithm also has a high degree of concurrency -- the low memory properties of the algorithm can be ``traded" for concurrency.
Exploiting these observations requires looking at the data dependencies of the FMG algorithm.

To fully exploit the available parallelism in the FMG algorithm we generalize the sweeping process of the original algorithm by defining a regular ``patch" $i$ of cells $u^k_i$, on grid $k$, with say 4-64 cells on a side.
Define a partitioning of each grid into a non-overlapping set of patches  $G^k$, and an extension of a patch $u_i$ by some number of cells as $\bar u_i$ -- these are \textit{halo or buffer cells} and they allow for subdomain solves with inaccurate boundary conditions to be solved accurately in the region of interest $u_i$ without communication.
These extended patches are conceptually similar to buffer regions that are used in algorithms to reduce the number of messages at the expense of sending more data and redundant work \cite{Leiserson:1997bc}.
Define a solver or smoother $S$ on an extended patch with non-homogenous boundary conditions that returns an improved solution on the original -- non-extended -- patch of data.
This smoother is used as the coarse grid solver for notational convenience and it must be accurate when used on the (entire) coarse grid. 
This smoother takes two extra arguments -- $I^l_{l+1} $ and ${\bar u}^{l}_i $ -- for use in a Kaczmarz smoother (\S \ref{ssec:Kaczmarz}).
The smoother assembles these patch or ``block" solve solutions \textit{additively}, in a block Jacobi method, to increase the degree of parallelism over the multiplicative method that is natural in serial.
%We also use the fourth order interpolation in the full update on the SR levels.
Assume that the coarsest grid, grid $u^0$, is composed of only one patch, again for notational convenience, and that each subsequent grid is a simple refinement by a small integer refinement ratio (ie, two or four).
The size of the group of patches on each level is a factor of eight times larger than the next coarser level in 3D in a non-AMR solve with a refinement ratio of two (or 64 with a refinement ratio of 4).
An AMR solve, with nested constant size patches, would pruned these groups appropriately.
Figure \ref{fig:fas-fmg-sr} shows a sketch of a \textit{parallel segmented refinement} algorithm assuming $M$ coarse grids and the forcing function $f$ has been suitably interpolated to, or defined on, all levels.

\vskip .1in
\begin{figure}[ht!]
\vbox{ \raggedright 
$\phantom{}${\bf  function} FMG-FAS-SR \\ 
\vskip .1in
$\phantom{MM}u^0 \leftarrow S\left(L^0, u^0, f^0\right)$ \\
$\phantom{MM}${\bf  for} $k=0:M-1$ \\
$\phantom{MMMM}${\bf for all} $u^k_i \in G^k$ \\ 
1:$\phantom{MMMMM}{\bar u}^{k+1}_i \leftarrow \Pi^{k+1}_k u^k$  \qquad \qquad \qquad  \qquad \qquad -- FMG prolongation \\
2:$\phantom{MMMMM}u^{k+1}_i \leftarrow S\left(L^{k+1},{\bar u}^{k+1}_i , {\bar f}^{k+1}_i\right)$ \\
3:$\phantom{MMMMM}{u}_{i}^k \leftarrow {\hat I}^k_{k+1} {\bar u}^{k+1}_i$ \qquad \qquad \qquad \qquad \qquad \qquad \quad -- restrict solution \\ 
4:$\phantom{MMMMM}{\tau}_{i}^k \leftarrow L^k \left( {\hat I}^k_{k+1} {\bar u}^{k+1}_i\right) - I^k_{k+1}\left(L^{k+1} {\bar u}^{k+1}_i\right)$ \quad -- data dependence for ${\bar u}^{k+1}$ \\
5:$\phantom{MMMMM}{f}^k_i  \leftarrow I^k_{k+1} {f}^{k+1}_i + {\tau}_{i}^k $ \\
$\phantom{MMMM}${\bf for all} $u^k_i \in G^k$ \\ 
$\phantom{MMMMMM}u^k_i \leftarrow S\left(L^k,{\bar u}^k_i , {\bar f}^k_i \right)$ \quad  \\
$\phantom{MMMM}${\bf  for} $l=k:-1:1$   \qquad \qquad \qquad \qquad \qquad \qquad \qquad -- pre-smoothing leg of V-cycle\\ 
$\phantom{MMMMMM}${\bf for all} $u^l_i \in G^l$ \\ 
$\phantom{MMMMMMMM}{u}_{i}^{l-1} \leftarrow {\hat I}^{l-1}_{l} {\bar u}^{l}_i$ \qquad \qquad \qquad \qquad \qquad \quad -- restrict solution \\ 
$\phantom{MMMMMMMM}{\tau}_{i}^{l-1} \leftarrow L^{l-1} \left( {\hat I}^{l-1}_{l} {\bar u}^{l}_i\right) - I^{l-1}_{l}\left(L^{l} {\bar u}^{l}_i\right)$ \\
$\phantom{MMMMMMMM}{f}^{l-1}_i  \leftarrow I^{l-1}_{l} {f}^{l}_i + {\tau}_{i}^{l-1} $ \\
$\phantom{MMMMMMMM}u^{l-1}_i \leftarrow S\left(L^{l-1},{\bar u}^{l-1}_i , {\bar f}^{l-1}_i \right)$ \quad  \\
$\phantom{MMMM}${\bf  for} $l=0:k$  \qquad \qquad \qquad \qquad  \qquad  \qquad  \qquad  \quad  -- post-smoothing leg of V-cycle \\ 
$\phantom{MMMMMM}${\bf for all} $u^l_i \in G^l$ \\
$\phantom{MMMMMMMM}$\textbf{if} $l=M-1$\\ 
$\phantom{MMMMMMMMMM}{\bar u}^{l+1}_i \leftarrow {\Pi}^{l+1}_l u^l$ \qquad \qquad   \quad \qquad -- using Eq. \ref{eq:fupdate2} \&  HO prolongation \\
$\phantom{MMMMMMMMMM}u^{l+1}_i \leftarrow S\left(L^{l+1},{\bar u}^{l+1}_i  , {\bar f}^{l+1}_i, I^l_{l+1}, {\bar u}^{l}_i  \right)$  \quad -- use CR \\
$\phantom{MMMMMMMMMM}$\textbf{compute functional of} $u^{l+1}_i$ \qquad  -- fine grid \\ 
$\phantom{MMMMMMMM}$\textbf{else} \\
$\phantom{MMMMMMMMMM}{\bar u}^{l+1}_i \leftarrow {\bar u}^{l+1}_i  + I^{l+1}_l \left( u^l - I_{l+1}^l u^{l+1} \right)$ \quad -- using Eq. \ref{eq:fupdate}  \\
$\phantom{MMMMMMMMMM}u^{l+1}_i \leftarrow S\left(L^{l+1},{\bar u}^{l+1}_i  , {\bar f}^{l+1}_i \right)$  \\
}
\caption{Segmented refinement FMG-FAS-SR algorithm (lines labeled 1-5 must be fused to avoid the need to store $u^{h}$)}
\label{fig:fas-fmg-sr}
\end{figure}

The dependency graph of this algorithm is similar to a forest of oct-trees with additional dependencies between neighboring trees.
If simple averaging is used in restriction then processing a coarse grid patch depends on the $R^D$ child patches (eg, an oct-tree), with a refinement ratio of $R$ in $D$ dimensions.
Higher order interpolation, which we use for prolongation, adds edges in the data dependency graph between the trees.

\begin{figure}[ht!]
\begin{center}
%\epsfile{file=./FMG.eps,scale=0.8}
\includegraphics[width=150mm]{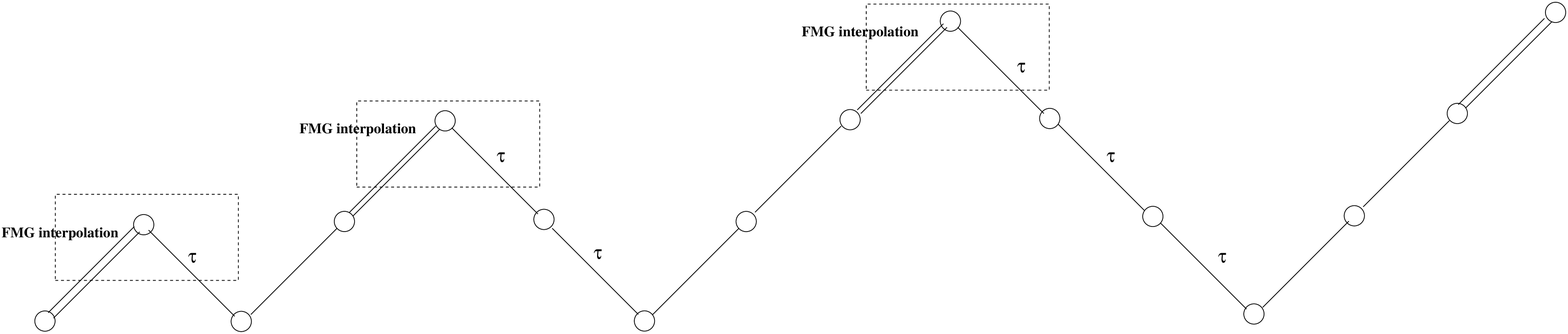} 
\caption{FMG cycle with $\tau$ corrections; dashed boxes show fused matrix free processing}
\label{fig:FMG}
\end{center}
\end{figure}

Figure \ref{fig:FMG} shows the FMG cycle, the $\tau$-corrections and the fine grid processing that can be fused and processed without permanent storage.
Note, that we use the higher order interpolation on the grids with full updates (SR levels), the finest level only in this figure.
This algorithm posses a high degree of concurrency, with, for instance, ten levels of refinement resulting in over one billion way parallelism in 3D and $R=2$.
%The algorithm does not become memory efficient until finer levels of refinement where a processing unit computes multiple patches sequentially.
%Generally one would not want to use SR on all levels, using normal FMG-FAS for the coarsest levels and using SR only on the finest grid or a few of the finest grids.

\subsection{Compatible Relaxation and Kaczmarz Smoother}
\label{ssec:Kaczmarz}

The critical change that we have made to the mathematical algorithm, to avoid storage of the finest grids, is the use of Eq. \ref{eq:fupdate2} to update the solution on the finest grids.
This method of not using a residual correction form, of using a full update, has the disadvantage that it adds coarse grid interpolation error to the whole fine grid solution instead of only to a correction.
%The post smoothing in the above algorithm will not maintain the correct solution at the coarse grid points if standard methods are used.
%This problem can be fixed by using a smoother that corrects for residuals of a certain equation that creep into the solution with standard smoothers.
We can ameliorate this problem by using stronger smoothers and using \textit{compatible relaxation} (CR).
CR uses a distributive relaxation or Kaczmarz relaxation in combination with a standard point-wise soother like Gauss-Seidel  \cite{Kaczmarz-37,Brandt-2011}.
Note, extra smoothing steps may be required, using extra flops, but because we have taken care to insure good data locality no additional memory movement is required, which is acceptable in the machine model that we are optimizing for.
We wish to maintain the approximation properties of the coarse grid while allowing smoothing of the error on the fine grid.
One approach for maintaining the approximation properties of the coarse grid is to (approximately) constrain the fine grid solution to solve
\begin{equation} 
{\tilde u}^H = I^H_h {\tilde u}^h.
\label{eq:kaz1}
\end{equation}
%\clearpage
Figure \ref{fig:smoother} shows the CR algorithm for our smoother on full update levels, that alternates between a standard smoother and a Kaczmarz relaxation.
\vskip .2in
\begin{figure}[ht!]
\vbox{ \raggedright 
$\phantom{}${\bf  function} $S(L^h,u^h,f^h,I^H_h,u^H)$ \\ 
$\phantom{MM} P \leftarrow I^H_h\left(I^H_h\right)^T $ \\ 
$\phantom{MM}${\bf for all} $j$ in patch $p^H$ \\ 
$\phantom{MMMM} r \leftarrow u^H - I^H_h(j,:)u^h$ \qquad  \quad  -- residual \\ 
$\phantom{MMMM} t \leftarrow r / P(j,j)$ \qquad  \qquad   \qquad  -- scalar correction \\ 
$\phantom{MMMM} u^h \leftarrow u^h + I_H^h(:,j)  t$ \qquad  \quad  -- update (distributive)\\ 
\vskip .1in
$\phantom{MM}${\bf standard smoother on patch} $u^h$ \\ 
}\caption{Kacmarz smoother on a patch}
\label{fig:smoother}
\end{figure}
%We use one Gauss-Seidel iteration as the standard smoother.
%\clearpage

%\clearpage
\subsection{The Solution}

A challenge of not explicitly storing the solution is the obvious problem of getting desired data from the simulation.
There are two basic methods for computing quantities of interest in the segmented refinement approach: 1) collect a functional of the data as the solution is computed, including streaming the entire solution to a ``file" for later processing (certainly useful for small simulations and debugging) and 2) storing a coarse grid solution which can be expanded or uncompressed efficiently with local processing (ie, PDE decompression) on demand for analysis.

%
%\bibliographystyle{siam}
%\bibliography{./the}
%\end{document}  

\section{Numerical Studies}
\label{sec:numr}

We investigate the properties of the algorithms developed here with a 1D Laplacian with homogenous Dirichlet boundary conditions and constant material coefficients.
We use a second order finite volume discretization, second order multigrid prolongation, fourth order FMG interpolation ($\Pi$) and first order accurate restriction operators.
The experiments are run in Matlab.
The Matlab source code is listed in \S\ref{sec:app}.
We do not use a data driven (vertical) processing that is a potential result of the algorithm in Fig. \ref{fig:fas-fmg-sr}, but our simple (horizontal) processing of the algorithm does have the same semantics as our proposed algorithm.
The smoother does simulate the asynchronous algorithm in that it is additive, a block Jacobi method, and so it is invariant to the order of processing of the blocks.
Each block has two non-overlapped cells, whose result is returned by the smoother, and two or four halo cells on each side (except at boundaries of course).
The solver within each subdomain is a few iterations of Gauss-Seidel, or compatible relaxation on the SR (full update) levels.
%We investigate normal FMG (no SR levels) and one, two, and three levels of SR.

To ascertain the costs of the proposed algorithm we conduct convergence studies, plotting the differential error $|\tilde{u} - u|_2$ as a function of the number of cells.
The discretization method is second order accurate and so we wish to maintain second order accuracy in our approximate solution with one FMG cycle.
We consider on 0-3 levels of SR in each study and look at the number of halo cells (two and four) and the number of Gauss-Seidel iterations (one and two) in the subdomain solver of the Jacobi smoother.
One application of the outer Jacobi smoother is used at all times.

\begin{figure}[ht!]
\begin{center}
\includegraphics[width=70mm]{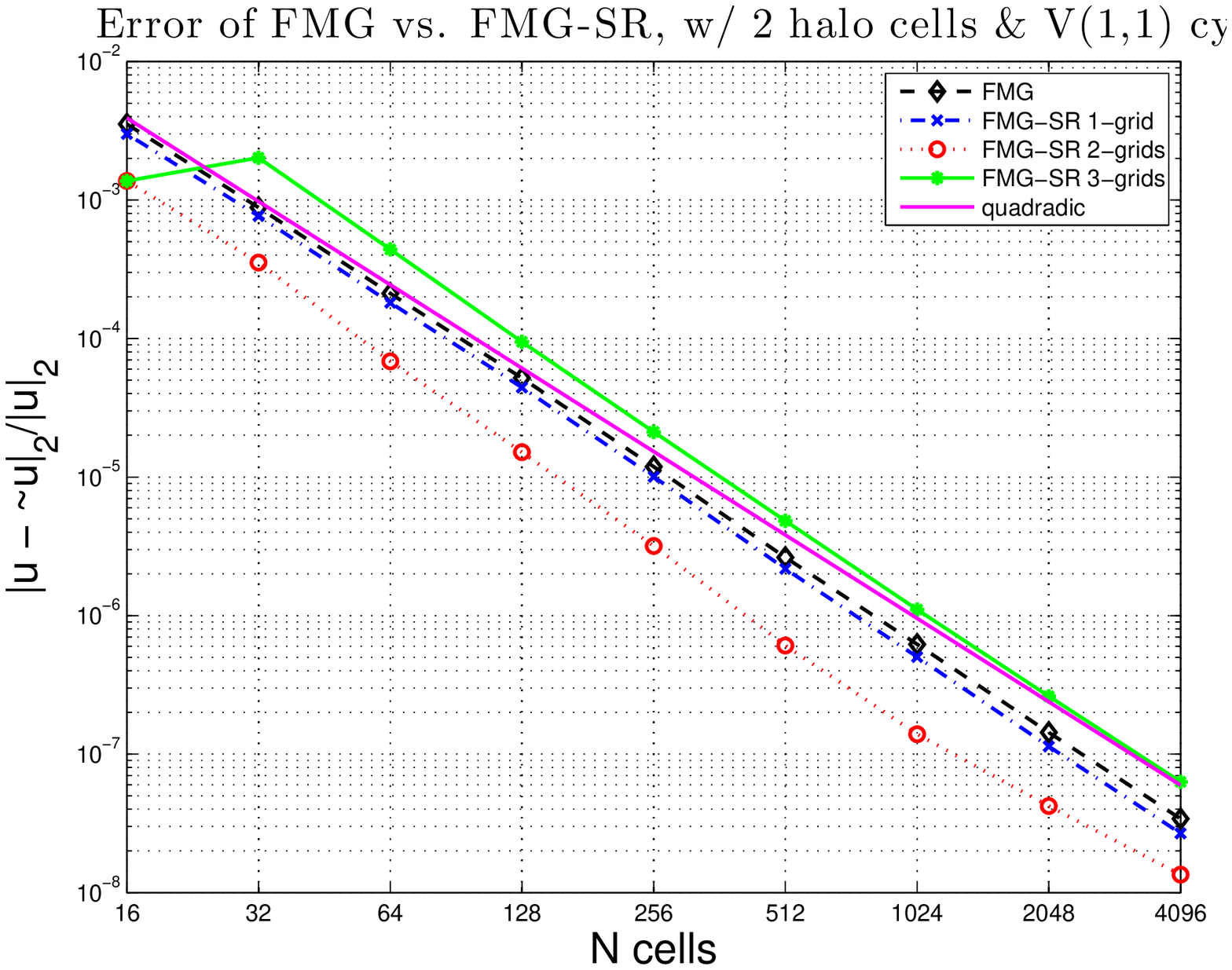} 
\includegraphics[width=70mm]{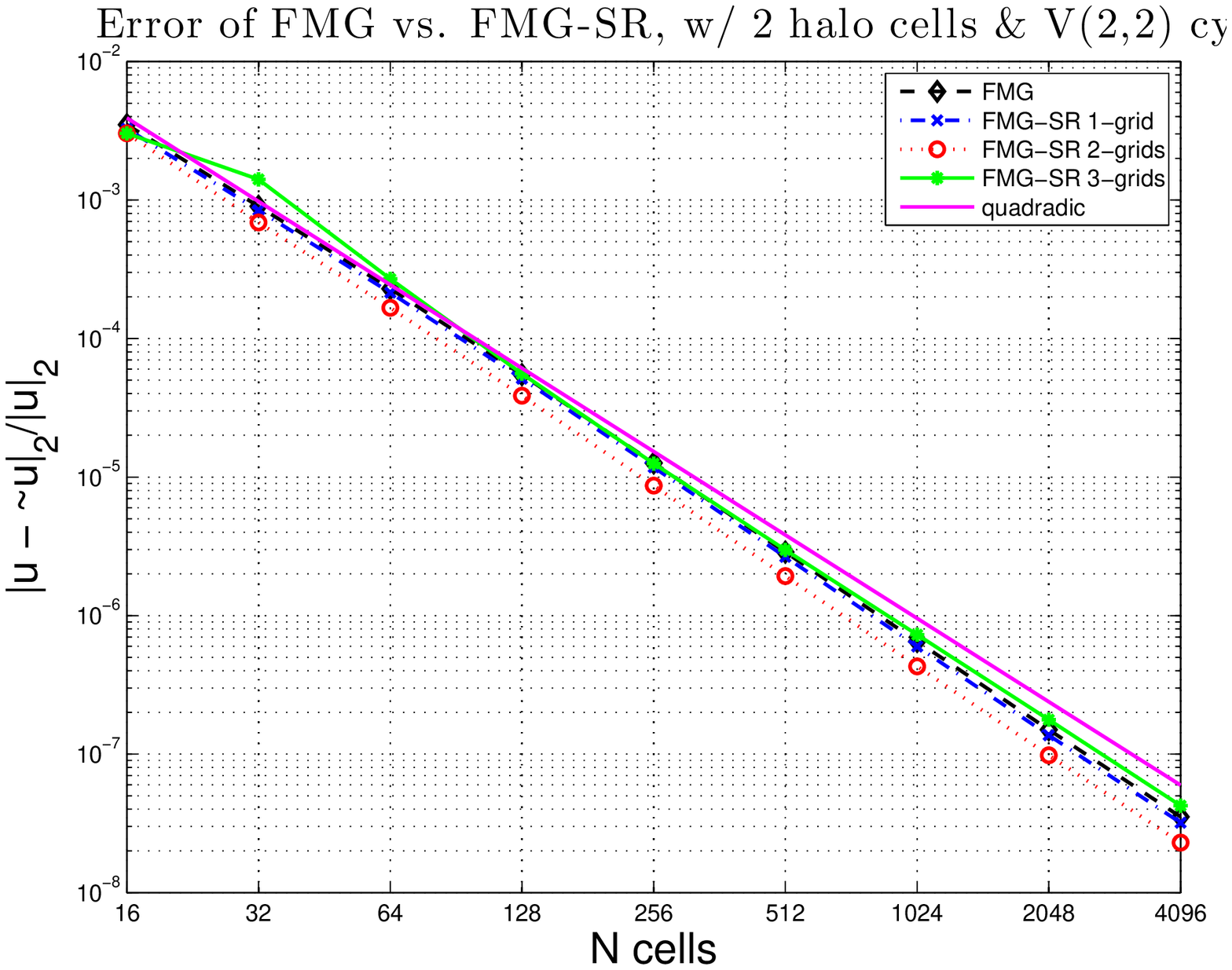} 
\caption{Convergence study with two halo cells in subdomains solver; V(1,1) cycle (left); V(2,2) cycle (right)}
\label{fig:conv_2halo}
\end{center}
\end{figure}

Figure \ref{fig:conv_2halo} shows convergence studies for the two halo smoother subdomains.
This data shows that truncation error accuracy (of the fine grid) is lost to some extent with three SR levels and $V(1,1)$ cycles, but otherwise we observe good second order convergence.

\begin{figure}[ht!]
\begin{center}
%\epsfile{file=./FMG.eps,scale=0.8}
\includegraphics[width=70mm]{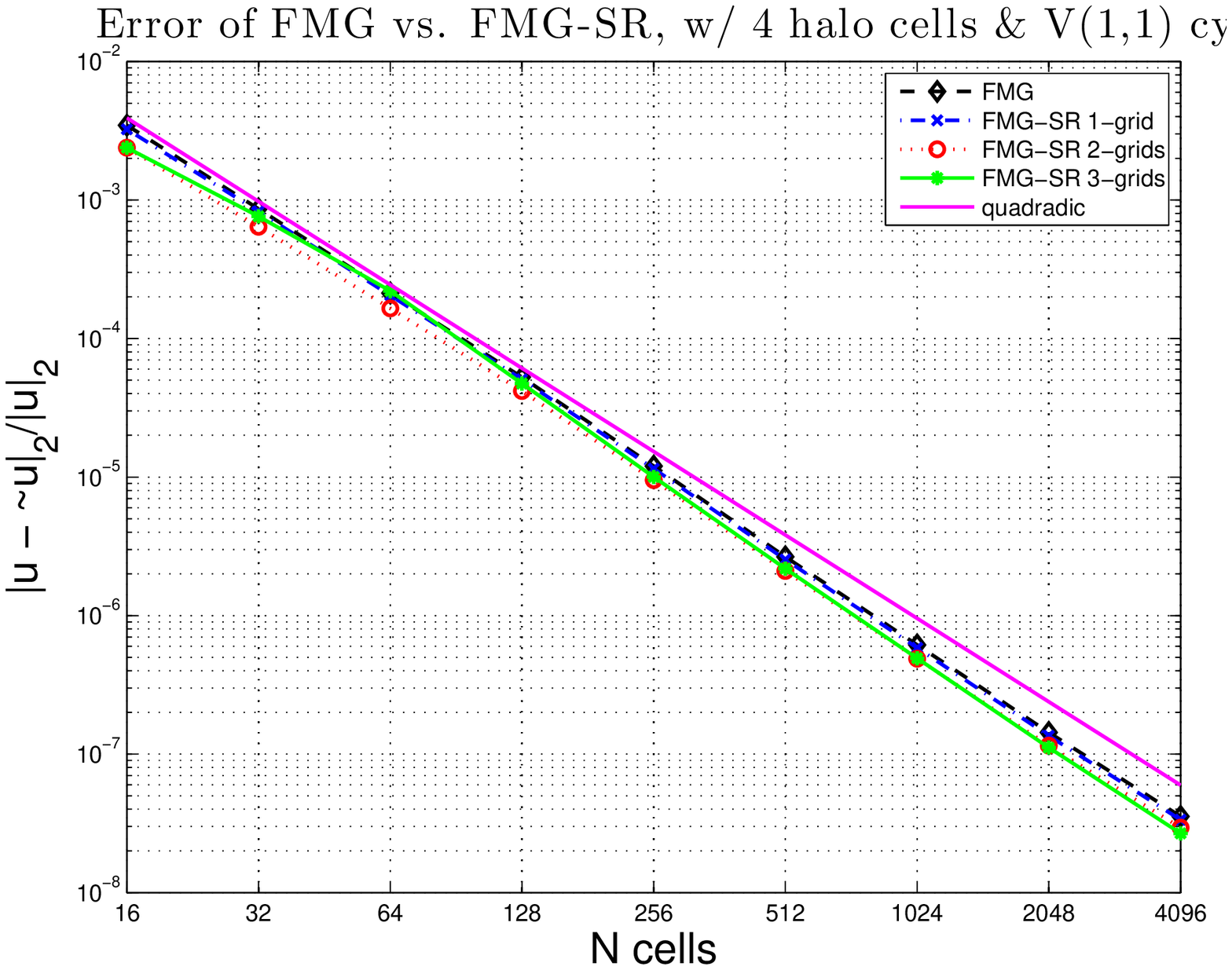} 
\includegraphics[width=70mm]{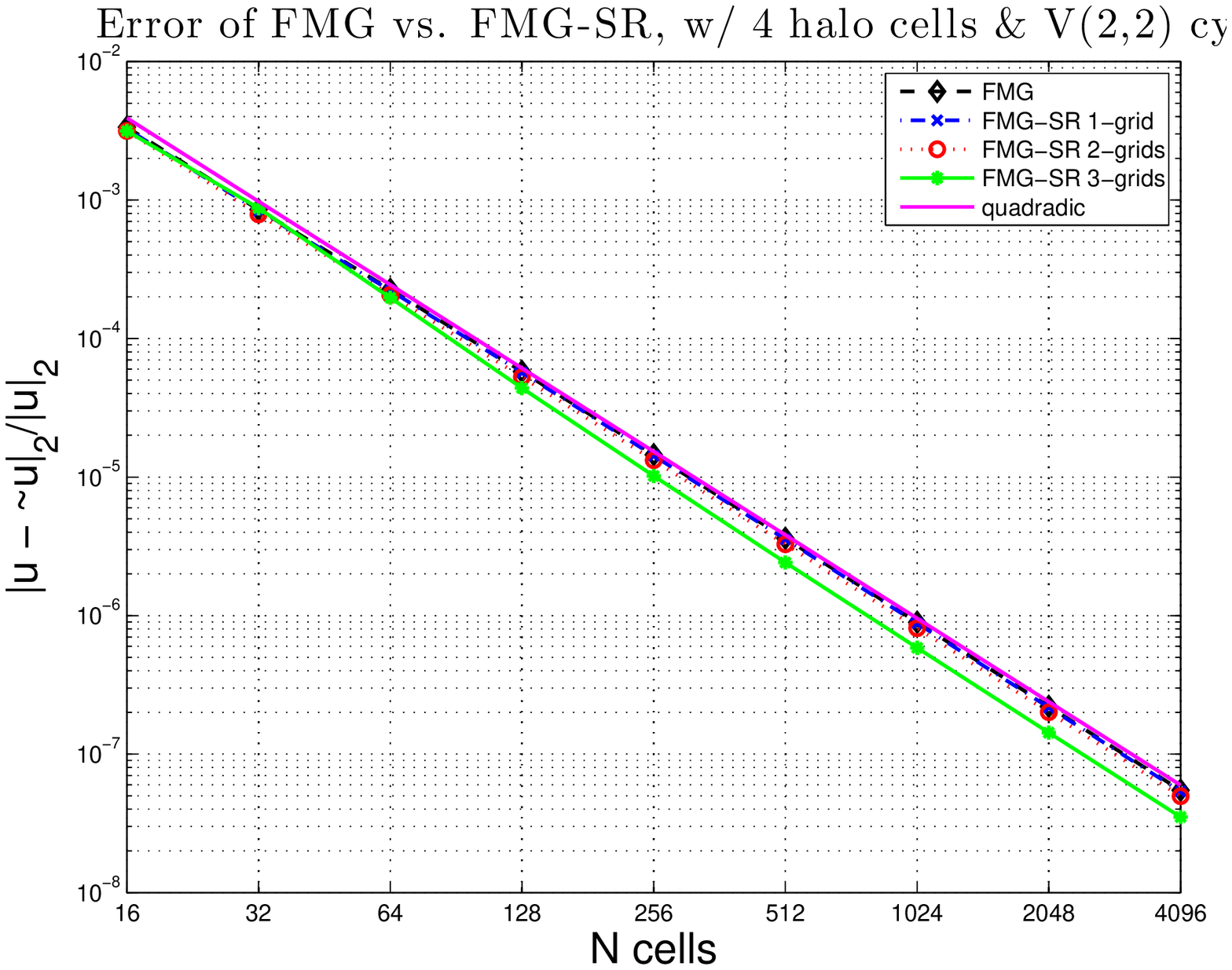} 
\caption{Convergence study with four halo cells in subdomains solver; V(1,1) cycle (left); V(2,2) cycle (right)}
\label{fig:conv_4halo}
\end{center}
\end{figure}

Figure \ref{fig:conv_4halo} shows convergence studies for the four halo smoother subdomains.
This data shows that with four halo cells the accuracy is very good with all solver configurations.

\begin{figure}[ht!]
\begin{center}
%\epsfile{file=./FMG.eps,scale=0.8}
\includegraphics[width=70mm]{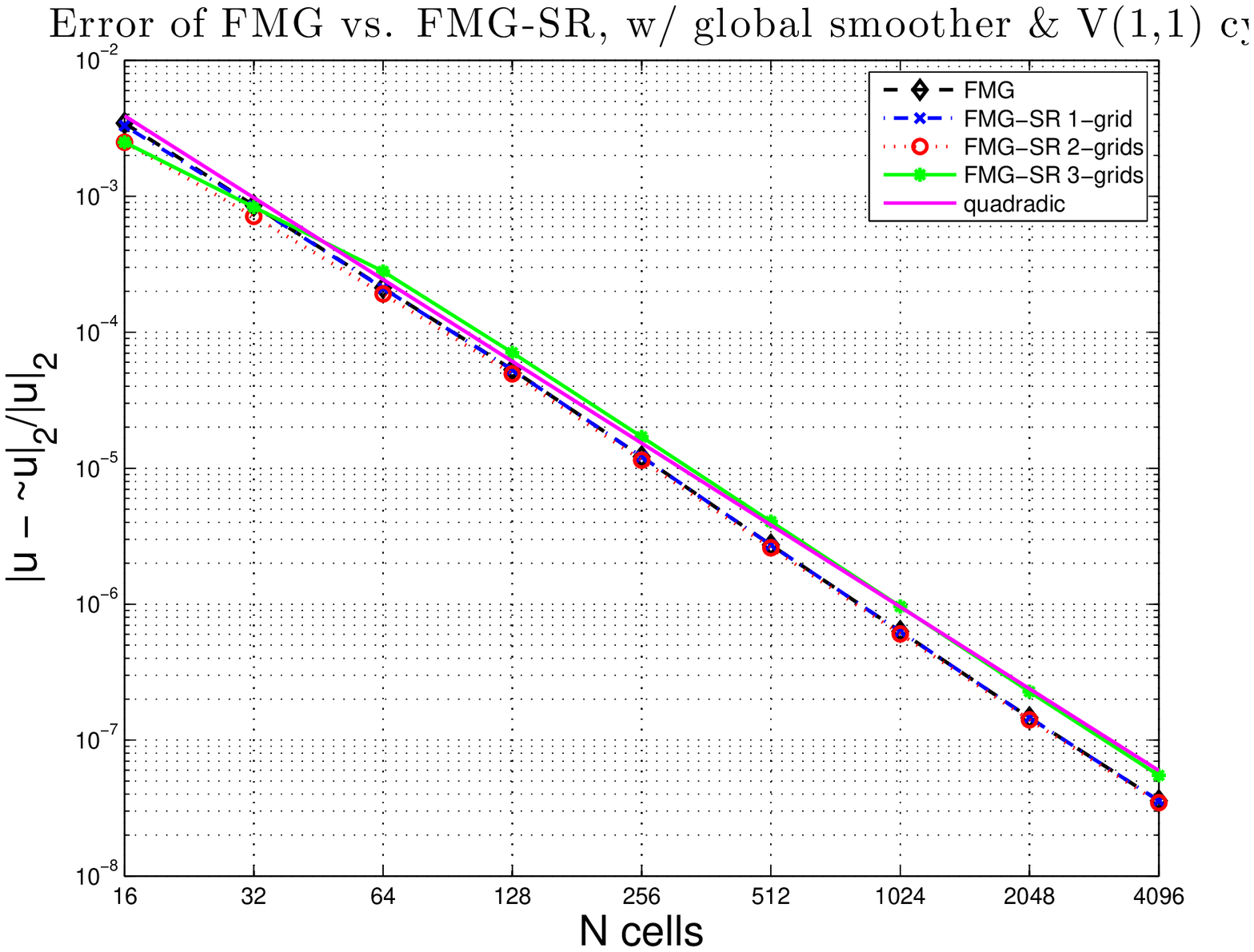} 
\includegraphics[width=70mm]{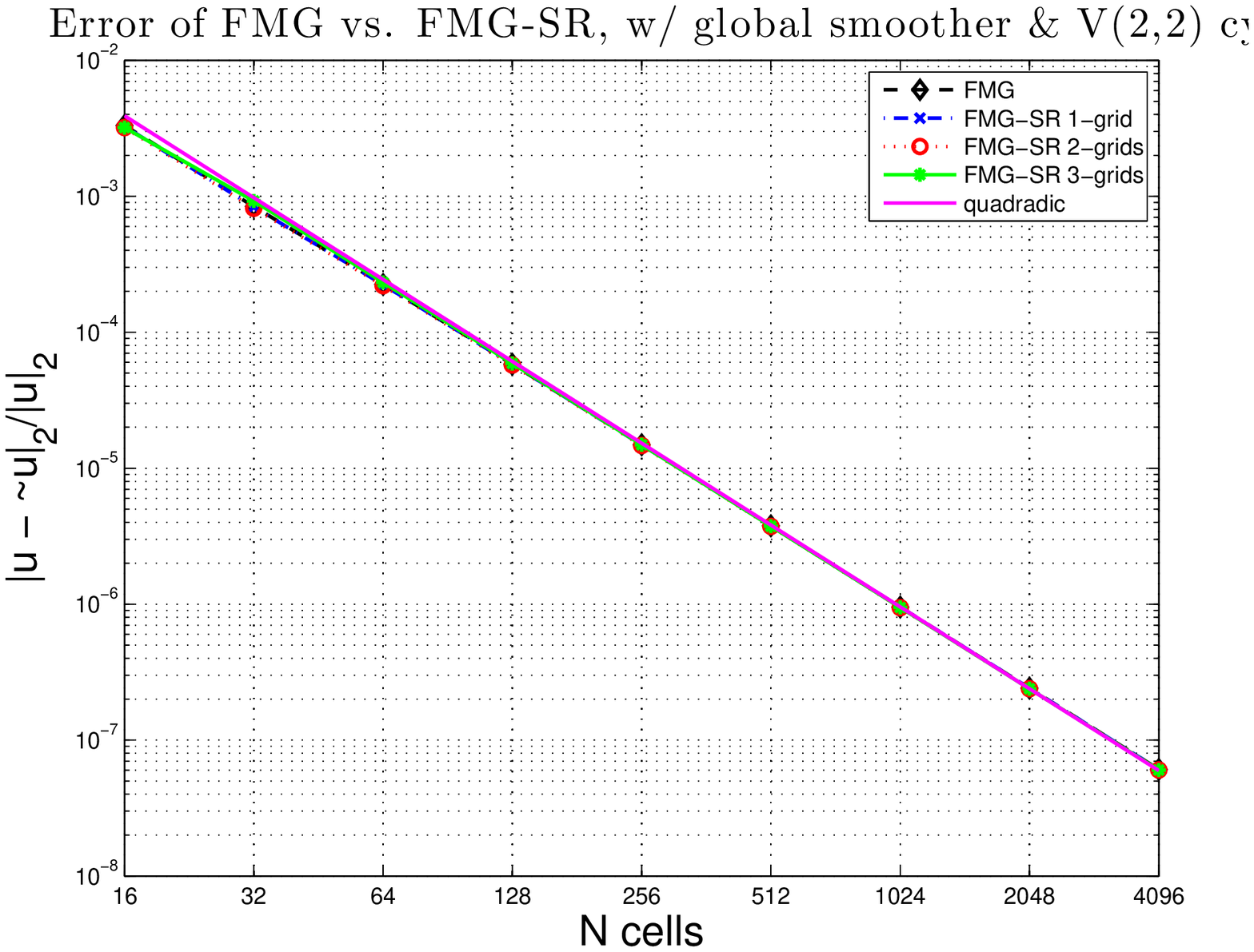} 
\caption{Convergence study with point-wise Gauss-Seidel smoothers; V(1,1) cycle (left); V(2,2) cycle (right)}
\label{fig:conv_globalhalo}
\end{center}
\end{figure}

For reference, Figure \ref{fig:conv_globalhalo} shows convergence studies using a simple point-wise Gauss-Seidel smoother (ie, the subdomain solver applied to the entire grid).
This data shows that the convergence results that we get, on this test problem, when using a standard (multiplicative) smoother are a bit ``cleaner" than that of the Jacobi smoothers in Figures \ref{fig:conv_2halo} and \ref{fig:conv_4halo}.

\section{Conclusions}
\label{sec:conc}

We have developed mathematical understanding of a highly concurrent FMG-FAS multigrid algorithm based on the $\tau$-correction and segmented refinement approach.
The method has the advantage of possessing very high levels of concurrency and is highly asynchronous.
This method also posses good data reuse properties because processing is confined to patches where operations can be ``fused", obviating the need to even store the entire solution at any one time.
We use overlapping subdomains which allows for accurate subdomain solves in the smoothers without communication.
These subdomain solves can be relatively accurate because the data is local (eg, in cache or fast memory of some sort) with low memory movement cost.
Interesting areas of future research are applying these methods to higher order discretizations, systems of PDEs, transient problems and hyperbolic problems, in a parallel and in an asynchronous environment.

\subsection*{Acknowledgments}
We would like to thank Achi Brandt for his generous guidance in developing these algorithms, Richard Vuduc for help in understanding advanced machine cost models, Hans Johansen for help with finite volume methods, and Jed Brown for many conversations on this topic.

\bibliographystyle{siam}
\bibliography{./low-mem-mg}

\section{Appendix}
\label{sec:app}
The data shown in this paper is generated with a Matlab code, shown here.
Run "conv(12)" to do a convergence study with 12 levels (4096 cells on the finest level).
The coarsest level is at level 2, so there are actually nine multigrid levels when ten levels are requested.

\begin{verbatim}
function conv( M )
% conv: Run convergence test for fmg()
%   Plot errors vs. N.
close all
set(0,'DefaultFigureWindowStyle','docked')
lw = 1.5; fz = 18;
s = 4; s_N = 2^s; s_h = 1/s_N;
nmodes = 16;
for halo_type=1:3
    for ns=1:2
        qq = s_h^2; 
        for k = s:M
            err_noRS(k)= fmg(k,0,2,0,ns,halo_type);
            err_RS1(k) = fmg(k,1,2,0,ns,halo_type);
            err_RS2(k) = fmg(k,2,2,0,ns,halo_type);
            err_RS3(k) = fmg(k,3,2,0,ns,halo_type);
            quad(k) = qq;
            qq = qq/4;
        end

        figure
        p=s:M;
        p2 = 2.^p;
        loglog(p2, err_noRS(4:M), 'kd--','linewidth',lw), hold on
        loglog(p2, err_RS1(4:M), 'bx-.','linewidth',lw), hold on
        loglog(p2, err_RS2(4:M), 'ro:','linewidth',lw), hold on
        loglog(p2, err_RS3(4:M), 'g*-','linewidth',lw), hold on
        loglog(p2, quad(4:M), 'm-','linewidth',lw), hold on
        set(gca,'XTick',p2)
        V = axis;
        V(1) = s_N; V(2) = 2^M;
        axis(V)
        legend('FMG','FMG-SR 1-grid','FMG-SR 2-grids','FMG-SR 3-grids','quadradic')
        ylabel('|u - ~u|_2/|u|_2','fontsize',fz);
        xlabel('N cells','fontsize',fz);
        if halo_type==3,
            title(['Error of FMG vs. FMG-SR, w/ global smoother \& V(',num2str(ns),',',num2str(ns),') cycles'],'fontsize',fz,'Interpreter','latex')
            grid
            print( gcf, '-djpeg100', ['conv_global_bc_',num2str(ns),'smooth'] )
            print( gcf, '-depsc', ['conv_global_bc_',num2str(ns),'smooth'] )
        else
            title(['Error of FMG vs. FMG-SR, w/ ',num2str(2*halo_type),' halo \& V(',num2str(ns),',',num2str(ns),') cycles'],'fontsize',fz,'Interpreter','latex')
            grid
            print( gcf, '-djpeg100', ['conv_',num2str(2*halo_type),'bc_',num2str(ns),'smooth'] )
            print( gcf, '-depsc', ['conv_',num2str(2*halo_type),'bc_',num2str(ns),'smooth'] )
        end
    end
end

function [ error ] = fmg( M, nRS, M0, pflag, ns, halo_type, nmodes )
% fmg: fmg-fas with segmetnted refinement.
%   M - number of coarse grids.
%   M0 - level of coarsest grid.
if nargin < 4, pflag = 1; end
if nargin < 3, M0 = 2; end
if nargin < 2, nRS = 0; end
if nargin < 1, M = 5; end
NN = 2^M, h_M = 1/NN;
if nargin < 7, nmodes = NN/16; end
if nargin < 5, ns = 1; end
if nargin < 6, halo_type = 2; end
set(0,'DefaultFigureWindowStyle','docked')

uu = cell(1,M);
rhs_orig = cell(1,M);
N = NN;
for k=M:-1:M0
 N_lev(k) = N; N = N/2;
end
%
[ L rhs ext Prol Rest RRt Prol_FMG ] = getops( M, nRS, M0, nmodes );
%
% FAS-FMG-SR w/ tau correction
%
prt = 0;
% FMG up to finest grid M
uu{M0} = smooth( L{M0}, 0, rhs{M0}, 1, M0 ); % coarsest grid solve
if prt, coarse_smoothing_error_infnorm = [M0 N_lev(M0) norm(uu{M0}-ext{M0},inf) norm(rhs{M0}-L{M0}*uu{M0},inf)], end
for k=M0:M-1
    uu{k+1} = Prol_FMG{k} * uu{k}; % FMG prol.
    % presmooth fines grid at this level
    if prt, pre_v_cycle_err_res_inf = [k+1 N_lev(k+1) norm(uu{k+1}-ext{k+1},inf) norm(rhs{k+1}-L{k+1}*uu{k+1},inf)], end
    rhs_orig{k+1} = rhs{k+1};
    uu{k+1} = smooth( L{k+1}, uu{k+1}, rhs{k+1}, ns, M0, 0, 0, 0, halo_type );
    % pre smoothing + coarse grid
    for m=k:-1:M0
        uu{m} = Rest{m}*uu{m+1}; % initial guess for coarse grid
        rhs{m} = Rest{m}*rhs{m+1} + L{m}*uu{m} - Rest{m}*L{m+1}*uu{m+1};
        uu{m} = smooth( L{m}, uu{m}, rhs{m}, ns, M0, 0, 0, 0, halo_type );
        %if prt, pre_smoothing_error_infnorm = [m N_lev(m) norm(uu{m}-ext{m},inf) norm(rhs{m}-L{m}*uu{m},inf)], end
    end
    % post smoothing
    for m=M0:k
        if m < M-nRS,
            uu{m+1} = uu{m+1} + Prol{m}*(uu{m} - Rest{m}*uu{m+1});
            uu{m+1} = smooth( L{m+1}, uu{m+1}, rhs{m+1}, ns, M0, 0, 0, 0, halo_type );
        else
            uu{m+1} = Prol_FMG{m}*uu{m};
            ns2 = ns; %[ 2*(M-m) N_lev(m+1) ], % (m-M+5)
            uu{m+1} = smooth( L{m+1}, uu{m+1}, rhs{m+1}, ns2(1), M0, uu{m}, diag(RRt{m}), Rest{m}, halo_type );
            %uu{m+1} = smooth( L{m+1}, uu{m+1}, rhs_orig{m+1}, ns2(1), M0, uu{m}, diag(RRt{m}), Rest{m} );
        end
        %if prt, post_smoothing_error_infnorm = [m+1 N_lev(m+1) norm(uu{m+1}-ext{m+1},inf) norm(rhs{m+1}-L{m+1}*uu{m+1},inf)], end
    end
    if prt, post_v_cycle_err_res_inf = [k+1 N_lev(k+1) norm(uu{k+1}-ext{k+1},inf) norm(rhs{k+1}-L{k+1}*uu{k+1},inf)], end
end
%res_red = norm(rhs{M}-L{M}*uu{M})/norm(rhs{M}),
%err_red = norm(uu{M}-ext{M})/norm(ext{M}),
%figure
%plot(rhs{M}-L{M}*uu{M},'b:*'), hold on,
%plot(rhs{M},'r:o'), hold on,
%pause
% plot & error
if pflag,
    close all
    figure
    xx = h_M/2 + h_M*(0:NN-1);
    plot(xx, uu{M}, 'r*--'), hold on
    plot(xx, rhs{M}, 'go--'), hold on
    plot(xx, ext{M}, 'bx-'), hold on
    plot(xx, abs(ext{M}-uu{M}), 'md-'), hold on
    axis([0 1 0 1.1*max(uu{M})])
    legend('result','b','x','error')
end
error = norm(ext{M}-uu{M},2)/norm(ext{M},2);

end

function [ L rhs ext Prol Rest RRt Prol_FMG ] = getops( M, nRS, M0, nmodes )
% getops: create opertors for FMG.
%   1D 2nd order finite volume discretization of Laplacian with Dirichlet
%   boundary conditions.  First order restriction, 2d oreder prolongation
%   and 4th order FMG interpolation.
NN = 2^M; h_M = 1/NN;

%
% Form restriction and prolongation ops
%
Prol = cell(1,M-1);
Rest = cell(1,M-1);
RRt = cell(1,M-1);
n=2^M0; m=2*n; 
for k=M0:M-1
  P = zeros(m,n); P0 = zeros(m,n);
  P(2,1) = 3; P(m-1,n) = 3;
  P(1,1) = 2; P(m,n) = 2;
  P0(2,1) = 1; P0(m-1,n) = 1;
  P0(1,1) = 1; P0(m,n) = 1;
  if m > 2,
      P(m-2,n) = 1; P(3,1) = 1; 
  end
  for j=2:n-1
      jj = (j-2)*2 + 2;
      pp = jj:jj+3;
      P(pp,j) = [ 1 3 3 1 ];
      pp = jj+1:jj+2;
      P0(pp,j) = [ 1 1 ];
  end
  Prol{k} = sparse(0.25*P);
  %Rest{k} = 0.125*P';
  Rest{k} = sparse(0.5*P0');
  RRt{k} = Rest{k}*Rest{k}';
  m = m*2; n = n*2;
end
%
% Form L & Prol_H
%
L = cell(1,M);
Prol_FMG = cell(1,M-1);
N = NN; h = h_M;
for k=M:-1:M0
 A = 2*eye(N) - diag(ones(N-1,1),1) - diag(ones(N-1,1),-1);
 A(1,1) = 3;
 A(N,N) = 3;
 L{k} = sparse(A*(1/h)^2);
 if k > M0,
    %Prol_H{k-1} = Prol{k-1};
    P = zeros(N,N/2);
    P(1,1) = 70;  P(1,2) = -2;
    P(2,1) = 112; P(2,2) = 35;   P(2,3) = -5;
    P(3,1) = 40;  P(3,2) = 105;  P(3,3) = -7;
    P(4,1) = -7;  P(4,2) = 105;  P(4,3) = 35;  P(4,4) = -5;
    %
    for i=5:2:N-4
        j = (i-1)/2 + 1;
        P(i,j-2) = -5;  P(i,j-1) = 35;   P(i,j) = 105;    P(i,j+1) = -7; 
                      P(i+1,j-1) = -7; P(i+1,j) = 105;  P(i+1,j+1) = 35;  P(i+1,j+2) = -5;
    end
    %
    j = N/2;
      P(N,j) = 70;    P(N,j-1) = -2; 
    P(N-1,j) = 112; P(N-1,j-1) = 35; P(N-1,j-2) = -5;
    P(N-2,j) = 40;  P(N-2,j-1) = 105; P(N-2,j-2) = -7;
    P(N-3,j) = -7;  P(N-3,j-1) = 105; P(N-3,j-2) = 35;  P(N-3,j-3) = -5;
    %
    Prol_FMG{k-1} = (1/128)*sparse(P);
 end
 N = N/2; h = h * 2;
end
%
% Form f
%
rhs = cell(1,M);
ext = cell(1,M);
N = NN; h = h_M;
x = h/2 + h*(0:N-1);
rhs{M} = 0*x'; ext{M} = 0*x';
for j=1:2:nmodes,
    rhs{M} = rhs{M} + (1/j)*sin(j*pi*x)';
    ext{M} = ext{M} + (1/j)*sin(j*pi*x)'/(j*pi)^2;
end
f = rhs{M};
e = ext{M};
%figure
%plot(x,(ext{M} - L{M}\rhs{M})./ext{M},'o--'), hold on
%plot(x,ext{M},'or--'), hold on
for k=M-1:-1:M0
    f = Rest{k} * f;
    rhs{k} = f;
    e = Rest{k} * e;
    ext{k} = e;
    %N = N/2; h = h * 2;
    %x = h/2 + h*(0:N-1);
    %plot(x,(ext{k} - L{k}\rhs{k})./ext{k},'*-.'), hold on
    %plot(x,ext{k},'*k-.'), hold on
end
%print( gcf, '-djpeg100', 'disc_error' ), grid, pause
end

function u_out = smooth( L, u_0, f, ns, M0, u_H, RRt, R, halo_type, a_p1, a_p2 )
% smooth: SR smoother
%
[N x] = size(L);
omega = 1.; 
sqrt2i = 1/sqrt(2);
if N == 2^M0, 
    % coarse grid solve
    u_out = L \ f;
else
    if nargin < 10,
        % whole level, additive Schwarz
        if halo_type==1,
            u_out = zeros(N, 1);
            p1 = 1:2;
            p2 = 1:4;
            u_out(p1) = smooth( L, u_0, f, ns, M0, u_H, RRt, R, halo_type, p1, p2 );
            for ii=6:2:N,
                p1 = (ii-3):(ii-2);
                p2 = (ii-5):ii;
                u_out(p1) = smooth( L, u_0, f, ns, M0, u_H, RRt, R, halo_type, p1, p2 );
            end
            p1 = (N-1):N;       
            p2 = (N-3):N;
            u_out(p1) = smooth( L, u_0, f, ns, M0, u_H, RRt, R, halo_type, p1, p2 );
        elseif halo_type==2,
            u_out = zeros(N, 1);
            p1 = 1:2;
            p2 = 1:6;
            u_out(p1) = smooth( L, u_0, f, ns, M0, u_H, RRt, R, halo_type, p1, p2 );
            p1 = 3:4;
            p2 = 1:8;
            u_out(p1) = smooth( L, u_0, f, ns, M0, u_H, RRt, R, halo_type, p1, p2 );
            for ii=10:2:N,
                p1 = (ii-5):(ii-4);
                p2 = (ii-9):ii;
                u_out(p1) = smooth( L, u_0, f, ns, M0, u_H, RRt, R, halo_type, p1, p2 );
            end
            p1 = (N-3):(N-2);
            p2 = (N-7):N;
            u_out(p1) = smooth( L, u_0, f, ns, M0, u_H, RRt, R, halo_type, p1, p2 );            
            p1 = (N-1):N;
            p2 = (N-5):N;
            u_out(p1) = smooth( L, u_0, f, ns, M0, u_H, RRt, R, halo_type, p1, p2 );
        else
            % whole level, G-S
            u_out = zeros(N, 1);
            p1 = 1:N; p2 = p1;
            u_out(p1) = smooth( L, u_0, f, ns, M0, u_H, RRt, R, halo_type, p1, p2 );            
        end
    else
        % one subdomain
          for k = 1:ns,
            % distributed (Kaczmarz) relaxation
            [n1 n2] = size(RRt);
            if n1+n2 > 2,
                % setup coarse grid iterator
                if a_p2(1) < a_p2(2), inc = 2; off = 0; else inc = -2; off = 1; end
                [x n] = size(a_p2);
                % distributed (Kacmarz) relaxation
                for ii=(a_p2(1):inc:a_p2(n)) - off
                    jj = (ii-1)/2 + 1;              % H index
                    r = u_H(jj) - R(jj,:)*u_0;
                    t = r / RRt(jj);             % update for RR' y = x^H
                    u_0 = u_0 + R(jj,:)' * t;       % update x^h
                end
                %error = [k size(a_p2) norm(u_H - R*u_0)]                
            end
            % normal smoothing
            for ii=a_p2
                u_0(ii) = u_0(ii) + omega * (f(ii) - L(ii,:)*u_0) / L(ii,ii);
            end
            % symmetrize
            a_p2 = fliplr(a_p2);
        end
        u_out = u_0(a_p1);
    end
end

end


\end{verbatim}

\end{document}